\documentclass[twoside, 11pt]{article}
\usepackage{mathrsfs}

\usepackage{mathrsfs,amsfonts,amsmath}
\RequirePackage{ifthen,calc}
\usepackage{theorem}
\usepackage[dvips]{color}

 \catcode`@=11
 \def\@evenhead{\hbox to\textwidth{\footnotesize\rm\thepage \hfill
  {\it }}} % authors name

 \def\@oddhead{\hbox to \textwidth{\footnotesize{\it
  Branching random walk in random environment with random absorption wall  } \hfill\thepage}}% abbreviate title

 \renewcommand{\section}{\makeatletter
 \renewcommand{\@seccntformat}[1]{{\csname the##1\endcsname.}\hspace{0.45em}}
 \makeatother \@startsection
{section}%                                            the name
{1}%                                                  the level
{0pt}%                                                the indent
{\baselineskip}%                                      the beforeskip
{0.5\baselineskip}%                                   the afterskip
{\normalsize\bfseries\mathversion{bold}}}
\catcode`@=12
\newcommand\ack{\section*{Acknowledgement}}

\newtheorem{thm}{\noindent Theorem}[section]
\newtheorem{lem}{\noindent Lemma}[section]
\newtheorem{cor}{\noindent Corollary}[section]

\theoremheaderfont{\normalfont\bfseries} \theorembodyfont{\slshape}
{\theorembodyfont{\rmfamily}}
{\theorembodyfont{\rmfamily}}
{\theorembodyfont{\rmfamily}}
{\theorembodyfont{\rmfamily}}
{\theorembodyfont{\rmfamily}}
{\theorembodyfont{\rmfamily}\newtheorem{rem}{\noindent Remark}[section]}
{\theorembodyfont{\rmfamily}}
{\theorembodyfont{\rmfamily}}
{\theorembodyfont{\rmfamily}}
 \def\beqlb{\begin{eqnarray}}\def\eeqlb{\end{eqnarray}}
 \def\beqnn{\begin{eqnarray*}}\def\eeqnn{\end{eqnarray*}}
 \numberwithin{equation}{section}
\def\L{{\mathcal L}}

\def\bfE{{\mathbb{E}}}
\def\bfP{{\mathbb{P}}}
\def\bfR{{\mathbb{R}}}

\def\bfN{{\mathbb{N}}}
\def\1{{\mathbf{1}}}

\begin{document}
\title{\bf Branching random walk in random environment with random absorption wall
}
\author{ You Lv\thanks{Email: youlv@mail.bnu.edu.cn }
\\ \small School of science of mathematics, Beijing Normal University,
\\ \small Beijing 100875, P. R. China.
}
\date{}
\maketitle

%\renewcommand{\thefootnote}{\fnsymbol{footnote}}\footnotetext[1]{}
%\renewcommand{\baselinestretch}{1.0}
%\noindent\hrulefill

\noindent\textbf{Abstract}: We consider the branching random walk in random environment with a random absorption wall. When we add this barrier, we discuss some topics related to the survival probability.  We assume that the random environment is i.i.d., $K_i$ is a particular i.i.d. random walk depend on the random environment $\L$. Let the random barrier function (the random absorption wall) is $g_i(\L):=ai^\alpha-K_i,$ where $i$ present the generation. We show that there exists a critical value $a_c>0$ such that if $a>a_c,\alpha=\frac{1}{3}$, the survival probability is positive almost surly and if $a<a_c,\alpha=\frac{1}{3},$
 the survival probability is zero almost surely. Moreover, if we denote $Z_n$ is the total populations in $n$-th generation in the new system (with barrier),
  under some conditions, we show $\ln\bfP_{\L}(Z_n>0)/n^{1/3}$ will converges to a negative constant almost surely if $\alpha\in[0,\frac{1}{3})$.
%We investigate the probability of a Brownian motion staying between two trajectories which are related to another Brownian motion.

\smallskip


\smallskip

\noindent\textbf{2000 Mathematics Subjects
Classification:} 60J80

%he environment is a sequence of point process laws
\section{Introduction}
The model named branching random walk on $\bfR$ with random environment in time (BRWre) has been introduced in \cite{HL2014} and \cite{MM2016}. Let $\L=(\L_1,\L_2,\ldots,\L_n,\ldots)$ be an i.i.d. random sequence of point process law which is also called the environment sequence. More precisely, $\L=(\L_1,\L_2,\ldots,\L_n,\ldots)$  is an i.i.d sequence of random
variables take values in the space of the distributions on the set of point processes on $\bfR$.
 After giving a realization $(L_1,L_2,\ldots,L_n,\ldots)$ of $\L,$ a time-inhomogeneous branching random walk is driven by the following way. It starts with
one individual located at the origin at time $0$.This individual dies at time 1 giving
birth to children and the children's position is according to the point process $L_1.$  Similarly, at each time $n$ every individual alive at generation $n-1$ dies and gives birth to children, and the position of the children with respect to their parent are given by the point process $L_n.$  We denote by $\mathbf{T}$ the (random) genealogical tree of the process. For a given individual $u\in\mathbf{T}$ we write $V(u)\in\bfR$ for the position of $u$ and $|u|$ for the generation at which $u$ is alive. The pair $(\mathbf{T},V)$ is called the branching random walk with i.i.d. random environment $\L.$  Conditionally on a realization environment sequence $\L$,  we denote $\bfP_{\L}$ for the law of this BRWre $(\mathbf{T},V)$ and $\bfE_{\L}$
for the corresponding expectation. The joint probability of the environment and the branching random walk is written $\bfP,$ with the
corresponding expectation $\bfE.$

Now we add an absorbing barrier to the BRWre. For a realization of environment $\L,$ we  write the barrier function $g_{\L}(i).$  At generation $i$, we erase all the individuals whose position is strictly larger than $g_{\L}(i)$ and its descendants.  We denote the new system $(\mathbf{T},V,g_{\L}(i)).$ we call it branching random walk with i.i.d. random environment and random absorbing barrier. This paper is focused on the survival or extinction problem when we add an absorbing barrier and the speed of extinction when we add a barrier which makes the system  $(\mathbf{T},V,g_{\L}(i))$ extinct.

When the environment space is degenerate, in other word, the branching random walk is time-homogeneous, the absorbing barrier problem has been researched by many scholars. Under the boundary case, Biggins et al \cite{BLSW1991} shows that if we let $g(i)=ai,$ then the system $(\mathbf{T},V,g(i))$ will survival if $a>0$ and extinct if $a\leq 0.$ Jaffuel \cite{BJ2012} gives a refinement order of critical barrier function, that is, if we let $g(i)=ai^\frac{1}{3},$ then there exist an $a_c>0$
 such that the system $(\mathbf{T},V,g(i))$ will survival if $a>a_c$ and extinct if $a< a_c.$ Furthermore, \cite{AJ2011} gives the speed of extinction when we take $g(i)\equiv 0$ and assume that the branching mechanism is $b~(b>2)$ binary tree.

For the model BRWre, Huang and Liu \cite{HL2014} proved that the maximal displacement in the process grows at ballistic speed
almost surely, and obtained central limit theorems and large deviations principles for the
counting measure of the process. Mallein \cite{MM2016} gives a more precise expression for the  asymptotic behaviour of maximal displacement which we will state in detail postponed.

\section{Basic assumption and main result}

First, we give some notations for the model BRWre. For every $n\in\bfN,$ let
$$\kappa_n(\theta):=\ln\bfE_{\L}\left(\sum_{l\in L_n}e^{-\theta l}\right),~~ \theta\in[0,+\infty)$$ be the log-Laplace transform of the point process $L_n,$ which is a point process according to the law $\L_n.$ We should notice that for fixed $\theta$ $\kappa_n(\theta)$ is also a random variable defined on environment space. Furthermore, $\{\kappa_n(\theta), n\in\bfN\}$ is an i.i.d. random sequence since the environment sequence is i.i.d.. We assume that $\mathbf{E}(\kappa^-_n(\theta))<+\infty$ for all $\theta\geq 0.$ Hence we can well define the function $\kappa: [0, +\infty)\rightarrow (-\infty,+\infty]$ by $\kappa(\theta):=\bfE (\kappa_n(\theta)).$ Then we can introduce four basic assumption in this paper.

{\it Condition 1: We assume that the interval~$\Upsilon:=\{\theta: |\kappa(\theta)|<\infty,~\kappa''(\theta)~\text{exists}\}$ has an non-empty inner and
 we can find a $\vartheta>0$ such that \beqlb\label{sect-2}\kappa(\vartheta)=\vartheta\kappa'(\vartheta).\eeqlb}
{\it Condition 2: There exists $\alpha_1,~\alpha_2,~\alpha_3>0$ such that
\beqlb\label{sect-4}\bfE(e^{\alpha_1|\vartheta \kappa'_1(\vartheta)-\kappa_1(\vartheta)|})<+\infty.\eeqlb
Denote $l_\vartheta:=l+\kappa'_1(\vartheta)$, we have $\bfP-\rm{a.s.}$
\beqlb\label{sect-5}
\frac{\alpha_2\bfE_{\L}\Big(\sum_{l\in L_1}|l_\vartheta|^3 e^{\alpha_2|l_\vartheta|-\vartheta l}\Big)}{\bfE_{\L}\Big(\sum_{l\in L_1} e^{-\vartheta l}\Big)}\leq \kappa''_1(\vartheta).\eeqlb
Moreover, we assume that
$$\bfP(\kappa''_1(\vartheta)>0)>0.$$}
There is also a condition 2 also has a more concise (but more stronger in fact) substitution. We write them as condition 2'.

{\it Condition 2': There exists $C\in\bfR, c,\lambda>0$ such that
\begin{eqnarray}\label{c2}
\max\{\kappa_1(\vartheta+\lambda\vartheta)-(1+\lambda)\kappa_1(\vartheta),~\kappa_1(\vartheta-\lambda\vartheta)-(1-\lambda)\kappa_1(\vartheta)\}<C,~{\rm \bfP-a.s.}.
\end{eqnarray} and
\begin{eqnarray}\label{c2a}
\kappa''_1(\vartheta)>c,  ~{\rm \bfP-a.s..}
\end{eqnarray}}

{\it Condition 3: Assume that $\kappa(0)>0$ and $\bfE\left(\sharp L_1\sum_{l\in L_1}e^{-\kappa_1(\vartheta)-\vartheta l}\right)<+\infty.$}
That is to say, we assume that the underlying branching process with random environment is supercritical, and the last condition is:

{\it Condition 4: There exists $x<0~,A\in\bfN$ such that  \beqlb\label{sect-6}\bfE\Big(\big[\ln \bfE_{\L}( 1_{\{\sum_{l\in L_1}\leq A\}} \sum_{l\in L_1}1_{\{\vartheta l+\kappa_1(\vartheta)\in[x,0]\}} )\big]^8\Big)<+\infty.\eeqlb}
By condition 2 or 2' we can know
\beqlb\label{sect-3}\sigma^2_Q:=\vartheta^2\bfE(\kappa''(\vartheta))\in(0,+\infty), \sigma^2_A:=\bfE((\kappa_1(\vartheta)-\vartheta\kappa'_1(\vartheta))^2)\in[0,+\infty).\eeqlb We can only see more information contained in condition 1-4 after introducing the many-to-onr formula in section 3.
%\alpha_2\vartheta \Big[e^{-\alpha_2(\kappa'_1(\vartheta)-\kappa_1(\vartheta))}\bfE_{\L}(\sum_{l\in L_1}l^3_\vartheta e^{2l_\vartheta})
%+e^{\alpha_2(\kappa'_1(\vartheta)-\kappa_1(\vartheta))}\bfE_{\L}(\sum_{l\in L_1}l^3_\vartheta)\Big].
We denote the function $\gamma$ is what we have defined in \cite[Theorem 2.1]{LY201801}. That is
$$\lim\limits_{t\rightarrow+\infty}\frac{-\ln P^0(\forall_{s\leq t} B_s\in[-\frac{1}{2}+\beta W_s, \frac{1}{2}+\beta W_s]|W)}{t}=\gamma(\beta),~~~\rm{a.s.}$$ where $B,W$ are two independent standard Brownian motions. From now on, we let $\beta:=\frac{\sigma_A}{\sigma_Q}$ and denote $\sigma_Q$ by $\sigma$ for simplicity.

\begin{thm}Define $$X_n(\L):=\sharp\{|u|=n:\forall i\leq n, V(u_i)\leq ai^{1/3}-\frac{K_i}{\vartheta}\},$$
$$\bfP_{\L, survive}=\bfP_{\L}(\exists u\in \mathcal{T}_{\infty},\forall i\geq 1, V(u_i)\leq ai^{1/3}-\frac{K_i}{\vartheta})$$
Under the condition 1-4, denote $K_i=\sum_{j=1}^{i}\kappa_j(\vartheta),$ The following statement is true.

a).When $a>\frac{3\sqrt[3]{6\gamma(\beta)\sigma^2}}{2\vartheta},$  $\bfP_{\L, survive}>0, \rm{\bfP-a.s.}.$ Moreover, the function $$\vartheta a=\vartheta b+\frac{3\gamma(\beta)\sigma^2}{b^2\vartheta^2}$$ has two solution $b_1,b_2$,  For any given $b\in(b_1,b_2),$ for any $\epsilon>0,$ there exist a large enough $M.$  we have
$$\bfP_{\L}\left(\varliminf\limits_{k\rightarrow+\infty}\frac{\ln X_{M^k}(\L)}{M^{\frac{k}{3}}}\geq b_2\vartheta-\epsilon\right)>0,~~\rm{\bfP-a.s.}.$$

b).When $a<\frac{3\sqrt[3]{6\gamma(\beta)\sigma^2}}{2\vartheta}$,  $\bfP_{\L, survive}=0.$
$\lim\limits_{n\rightarrow\infty}\frac{\ln\bfP_{\L}(X_n(\L)>0)}{\sqrt[3]{n}}=c,~\rm{\bfP-a.s.},$ where $c$ is a negative constant dependent on $a$.
\end{thm}

\begin{thm}
Define $g(\cdot):\bfN\rightarrow\bfR^+$ satisfied that $\lim\limits_{n\rightarrow\infty}\frac{\sup_{i\leq n}g(n)}{n^{1/3}}=0,$ $X_n(\L)$ is the surviving population of the generation $n$ in the system $(\mathbf{T},V,g_{\L}(i)).$  That is to say, $$X_n(\L):=\sharp\{|u|=n:\forall i\leq n, V(u_i)\leq g(i)-\frac{K_i}{\vartheta}\}.$$
Under the condition 1-4, we have $\lim\limits_{n\rightarrow\infty}\frac{\ln\bfP_{\L}(X_n(\L)>0)}{\sqrt[3]{n}}=-\sqrt[3]{3\sigma^2\gamma(\beta)}.~\rm{\bfP-a.s.}.$
\end{thm}
\begin{rem}
Under some assumptions, according the result of \cite{MM2016}, we have $$\lim\limits_{n\rightarrow+\infty}\frac{\min_{|u|=n}V(x)+\frac{K_n}{\vartheta}}{\ln n}=c,~~~~~\text{in Probability}$$
We can see in the random environment, the first order of asymptotic behavior of the leftmost position is a random walk, that is why we set the
barrier function like $g_{\L}(i)=g(i)-\frac{K_i}{\vartheta}.$
\end{rem}
%\noindent{\bf Proof of  3.1}

%c stand for two positive constants respectively large enough and small enough, that may change from line to line.

\section{Some useful lemma}%Small deviation estimation in random environment
Now we introduce some useful lemmas.

\noindent{\bf Bivariate version many-to-one formula in random environment}

Many to one formula is essential in studies of extremal behaviour of branching random walks. The random environment version of Many to one formula has been first introduced in \cite{BM2015}. When the environment is degenerate, the bivariate version many-to-one formula can be found in \cite{GHS2011}.  In this paper we need a bivariate version many-to-one formula in random environment.  For every $n\geq 1,$ we write $L_n$ for a realisation of the point process with law $\L_n$.  Let $(X_i,\xi_i)$ be a random variable taking values in $\bfR\times\bfN$ such that for any measurable nonnegative function $f,$
\beqlb\label{m1}\bfP_{\L}(X_i\leq x, \xi_i\leq A)=\bfE_{\L}(1_{\{\sum_{l\in L_i}1\leq A\}}\sum_{l\in L_i}1_{\{l\leq x\}}e^{-\theta l-\kappa_i(\theta)}).\eeqlb
%$$\kappa_{n}(\vartheta)=\ln\bfE_{\L}(\sum_{l\in L_n}e^{-\vartheta l})$$
%$$\mu_n\big((-\infty,x]\big)=\bfE_{\L}(\sum_{l\in L_n}\1_{\{l\leq x\}}e^{-\vartheta l-\kappa_{n}(\vartheta)})=\frac{\bfE_{\L}(\sum_{l\in L_n}\1_{\{l\leq x\}}e^{-\vartheta l})}{\bfE_{\L}(\sum_{l\in L_n}e^{-\vartheta l})}, ~~~ \forall x\in\bfR.$$
So in quenched sense, $\{X_n\}_{n\in\bfN}$ is a sequence of independent random variables. We set $S_n=S_0+\sum_{i=1}^{n}X_i, S_0=0.$
%From now on, $\bfP_{\L}$ stands for the joint law of the BRWre ($\mathbf{T},V$) and the random variable $(X_i, \xi_i),$ conditionally on the environment $\L.$
For the inhomogeneous, we need to introduce the shift operator~$\mathfrak{T},$~define
$$\mathfrak{T}\L:=(\L_2,\L_3,\ldots), ~~\mathfrak{T}_k:=\mathfrak{T}^{*k},~~\mathfrak{T}_0\L:=\L.$$
That is to say~$\mathfrak{T}_k\L=(\L_{k+1},\L_{k+2},\ldots).$ we use~$\bfP^k_{\L}$ to present the distribution of~$(\mathcal{T},V,\bfP_{\mathfrak{T}_k\L})$ The corresponding expectation of~$\bfP^k_{\L}$ is~$\bfE^k_{\L}.$ For the consistency of the notation, we agree $S_n=\sum_{i=1}^{n}X_{k+i}$~ under $\bfP^k_{\L}.$% 而概率~$\bfP^k_{\L}$ 下的
Writing~$\{\xi_{n}\}_{n\in\bfN}$ to present $\{\xi_{k+n}\}_{n\in\bfN}$ under~$\bfP_{\L}.$ Hence we can write \eqref{m1} as
\beqlb\bfP_{\L}(X_{i+1}\leq x, \xi_{i+1}\leq A)&=&\bfP^i_{\L}(X_{1}\leq x, \xi_{1}\leq A)\nonumber\\&=&\bfE^i_{\L}(1_{\{\sum_{l\in L_1}1\leq A\}}\sum_{l\in L_1}1_{\{l\leq x\}}e^{-\theta l-\kappa_{i+1}(\theta)}).\eeqlb
The following fomula give the relationship between BRWre and RWre.
\begin{lem}
$(${\bf many-to-one}$)$~For any $k,n,A_i\in\bfN, i\leq n$ and any measurable non-negative function $f:(\bfR^n\times\bfN^n)\rightarrow \bfR,$ we have
\beqnn
&&\bfE^k_{\L}\left[\sum_{|u|=n}f(V(u_i),1\leq i\leq n)\1_{\{\gamma(u_{i-1})\leq A_i,1\leq i\leq n\}}\right]\nonumber
\\&&~~~~~~~~~~~~~~~~=\bfE^k_{\L}\left[e^{\vartheta S_n+\sum^{n}_{i=1}\kappa_{k+i}(\vartheta)}f(S_i,1\leq i\leq n)\1_{\{\xi_{i}\leq A_i,1\leq i\leq n\}}\right].~\rm{\bfP-a.s.}
\eeqnn
where $\gamma(u)$ present the children number of $u.$
\end{lem}
The proof of this lemma can be done by induction on $n,$ which is standard, just like the proof of \cite[Theorem 1.1]{S2015}. here we omit it.

{\bf Mogul'ski\v{\i} estimation}
Mogul'ski\v{\i} estimation is also an essential tool in the barrier problem of Branching random walk. Mogul'ski\v{\i} estimation had first introduced in \cite{Mog1974}. Mallein \cite{Mal2015} gives the time inhomogeneous version of Mogul'ski\v{\i} estimation. Here we give the random environment version of  Mogul'ski\v{\i} estimation (Lemma 3.2), the proof of Lemma 3.2 can be found in \cite{Lv201802}.  Let $T_n$ is a random walk with i.i.d. random environment in time and satisfied the following assumption. We use $\mu=(\mu_1,\mu_2,\cdots,\mu_n,\cdots)$ to present the random environment. Denote $$M_{n}:=\bfE_\mu(T_n), ~U_{n}:=T_n-\bfE_\mu(T_n),~\Gamma_n:=\bfE_\mu(U^2_n)=\bfE_\mu(T^{2}_n)-M^2_{n}.$$
 \begin{itemize}
  \item[(H1)] $\bfE M_1=0,\sigma^2_A:=\bfE(M_1^2)\in[0,+\infty),\sigma^2_Q:=\bfE(U_1^2)=\bfE(\Gamma_1)\in(0,+\infty).$
  \item[(H2)] There exists $\lambda_1>0,$ such that $\bfE(e^{\lambda_1 |M_{1}}|)<+\infty.$
  \item[(H3)] There exist $\lambda_2, \lambda_3>0$ such that $\bfE_\mu(e^{\lambda_2|U_1|})\leq \lambda_3\leq\bfE_\mu(U^2_1)$ almost surely.
  \item[(H4)] $\{r_n\}_{n\in\bfN}$ is a positive sequence such that for any $\varrho>0, \lim\limits_{n\rightarrow+\infty}\frac{n^{\varrho}}{r_n}=0.$ $f(n)$ is an positive integer-valued function such that for any $\kappa>0, \lim\limits_{n\rightarrow+\infty}\frac{f(n)}{e^{n^{\kappa}}}=0.$
\end{itemize}
Let $\xi_i$ be a positive random variable whose law is only determined by the $i-$th element $\mu_i$ in a realistic of environment $\mu.$ Moreover,
conditioned on a given environment realization $\mu,$  $\{\xi_i\}_{i\in\bfN}$ is an independent positive random sequence. We can also know for any measurable function $\eta$, $\{\bfE_{\mu}(\eta(\xi_i))\}_{i\in\bfN}$ is an i.i.d. random sequence in the environment space since $\mu$ is i.i.d..
\begin{lem}
Under the assumption $\bfE(\xi_1)<+\infty$ and (H1)-(H4), let $g(s), h(s)$ be two continue functions on $[0,1]$ and $g(s)<h(s)$ for any $s\in [0,1].$ $g(0)< a_0\leq b_0 <h(0), g(1)\leq a'<b'\leq h(1).$
Denote $C^{z_1,z_2}_{g,h}:=\int_{z_1}^{z_2}\frac{1}{[h(s)-g(s)]^2}ds,$ then for any $\alpha\in (0,\frac{1}{2}),$ then ${\rm \bfP-a.s.}$ we have
$$\varlimsup\limits_{n\rightarrow +\infty}\sup\limits_{x\in\bfR}\frac{\ln\bfP_\mu
(\forall_{0\leq i\leq n}\frac{T_{f(n)+i}}{n^\alpha}\in [g(\frac{i}{n}),h(\frac{i}{n})]|T_{f(n)}=x)}{n^{1-2\alpha}}\leq -C^{0,1}_{g,h}\sigma^2_Q\gamma(\frac{\sigma_A}{\sigma_Q}),$$
\beqnn&&\varliminf\limits_{n\rightarrow +\infty}\inf\limits_{x\in[a_0 n^{\alpha}, b_0 n^{\alpha}]}\frac{\ln \bfP_\mu
\Big(\substack{\forall_{0\leq i\leq n},T_{f(n)+i}\in [g(\frac{i}{n})n^\alpha,h(\frac{i}{n})n^\alpha]\\T_{f(n)+n}\in [a'n^\alpha,b'n^\alpha],~\xi_{i+f(n)}\leq r_n}\Big|T_{f(n)}=x\Big)}{n^{1-2\alpha}}
\\&&~~~~~~~~~~~~~~~~~~~~~~~~~~~~~~~~~~~~~~~~~~~~~~~~~~~~~~~~~~~~~~~~~~~~~~\geq -C^{0,1}_{g,h}\sigma^2_Q\gamma(\frac{\sigma_A}{\sigma_Q}).
\eeqnn
\end{lem}
\begin{cor}
Under the assumption of Lemma 3.2, $0\leq l< m\leq N.$ ${\rm \bfP-a.s.}$ we have
$$\varlimsup\limits_{k\rightarrow +\infty}\sup\limits_{x\in\bfR}\frac{\ln\bfP_\mu
(\forall_{lk\leq i\leq mk}\frac{T_{i}}{(Nk)^\alpha}\in [g(\frac{i}{Nk}),h(\frac{i}{Nk})]|T_{lk}=x)}{(Nk)^{1-2\alpha}}\leq -C^{\frac{l}{N},\frac{m}{N}}_{~g,~h}\sigma^2_Q\gamma(\frac{\sigma_A}{\sigma_Q}),$$
\beqnn&&\varliminf \limits_{n\rightarrow +\infty}\inf\limits_{x\in[a_0 n^{\alpha}, b_0 n^{\alpha}]}\frac{\ln \bfP_\mu
\left(\forall_{lk\leq i\leq mk}\substack{\frac{T_{i}}{(Nk)^\alpha}\in [g(\frac{i}{Nk}),h(\frac{i}{Nk})],\\\frac{T_{mk}}{(Nk)^\alpha}\in [a',b'],\xi_i\leq r_{Nk}}|T_{lk}=x\right)}{(Nk)^{1-2\alpha}}
\\&&~~~~~~~~~~~~~~~~~~~~~~~~~~~~~~~~~~~~~~~~~~~~~~~~~~~~~~~\geq -C^{\frac{l}{N},\frac{m}{N}}_{~g,~h}\sigma^2_Q\gamma(\frac{\sigma_A}{\sigma_Q}).
\eeqnn
\end{cor}
{\bf Proof of Corollary 3.1}

When $l=0,m=N,$ it has contained in Lemma 3.2.  Let $n:=mk-lk,  f(n)=\frac{ln}{m-l}=lk.$ we note that $\frac{T_{i}}{(Nk)^\alpha}\in [g(\frac{i}{Nk}),h(\frac{i}{Nk})]$ is equal to
{\small $$\frac{T_i}{(mk-lk)^\alpha}\in\left[\left(\frac{Nk}{mk-lk}\right)^\alpha g\Big(\frac{i-lk+lk}{mk-lk}\frac{mk-lk}{Nk}\Big),\left(\frac{Nk}{mk-lk}\right)^\alpha h\Big(\frac{i-lk+lk}{mk-lk}\frac{mk-lk}{Nk}\Big)\right].$$}
We still let $\beta:=\frac{\sigma_A}{\sigma_Q}$ and denote $\sigma_Q$ by $\sigma$ for simplicity. Let $X:=(x+\frac{lk}{mk-lk})\frac{mk-lk}{Nk},$ by lemma 3.2 we can see
\beqnn&&\varlimsup\limits_{k\rightarrow +\infty}\frac{\sup\limits_{x\in\bfR}\ln\bfP_\mu
(\forall_{lk\leq i\leq mk}\frac{T_{i}}{(Nk)^\alpha}\in [g(\frac{i}{Nk}),h(\frac{i}{Nk})]|T_{lk}=x)}{(mk-lk)^{1-2\alpha}}\frac{(mk-lk)^{1-2\alpha}}{(Nk)^{1-2\alpha}}
\\&=&\gamma(\beta)\sigma^2
\int_{0}^{1}\Big[(\frac{mk-lk}{Nk})^{\alpha}(h(X)-g(X))\Big]^{-2}dx\times\frac{(mk-lk)^{1-2\alpha}}{(Nk)^{1-2\alpha}}
\\&=&\gamma(\beta)\sigma^2\int_{0}^{1}\Big[h(X)-g(X)\Big]^{-2}dX
\\&=&\gamma(\beta)\sigma^2\int_{\frac{l}{N}}^{\frac{m}{N}}(h(x)-g(x))^{-2}\eeqnn
%$$H(x):=h\Big((x+\frac{lk}{mk-lk})\frac{mk-lk}{Nk}\Big), $$
The following lemma is considering the case $h(0)=b_0.$
\begin{cor}

\

\begin{itemize}
  \item[(H1)]Let $\nu\in(\alpha,1),$ $h(t)$ satisfied that  $\varliminf\limits_{t\rightarrow 0}\frac{h(t)-h(0)}{t^{\nu}}>-\infty.$
There exist a pair of $x<y<0$ such that $$\bfE([-\ln\bfP_\mu(T_{1}\in [x,y],\xi_{1}\leq A|T_{0}=0)]^8)<+\infty.$$
  \item[(H2)]$\varliminf\limits_{t\rightarrow 0}\frac{h(t)-h(0)}{t^{\alpha}}:=d>-\infty.$ there exists $x<d$ such that
  $$\bfE([-\ln\bfP_\mu(T_{1}\in [x,d),\xi_{1}\leq A|T_{0}=0)]^8)<+\infty.$$
  \end{itemize}
  If assumption (H1) or (H2) holds, $\alpha\leq 1/3,$ $a'<b'.$ we have
\beqnn\varliminf \limits_{n\rightarrow +\infty}\frac{\ln \bfP_\mu
\Big(\substack{\forall_{f(n)\leq i\leq f(n)+n}\frac{T_{i}}{n^\alpha}\in [g(\frac{i-f(n)}{n}),h(\frac{i-f(n)}{n})],\\ \frac{T_{f(n)+n}}{n^\alpha}\in [a',b'],\xi_i\leq r_{n}|T_{f(n)}=h(0)n^{\alpha}}\Big)}{n^{1-2\alpha}}
\geq -C^{0,1}_{g,h}\sigma^2_Q\gamma(\frac{\sigma_A}{\sigma_Q}).\eeqnn
\end{cor}
{\bf Proof of Corollary 3.2} Without loss of generality, we assume that $H1$ holds. By $\varliminf\limits_{t\rightarrow 0}\frac{h(t)-h(0)}{t^{\nu}}>-\infty, $
 we have $h(t)-h(0)>-dt^{\nu}$ when $t$ is small enough.  Choose a $\delta>0$ arbitrarily such that $\frac{g(0)-h(0)}{y}>\delta.$~(Under assumption (H2),we should need $\frac{g(0)-h(0)}{d}>\delta.$) Let $N=\lfloor\delta n^{\alpha}\rfloor.$ For the continuity of $g$, we can see for any small enough $\epsilon>0$ such that $g(0)+\epsilon-h(0)\leq y\delta$ then we can find a large enough $n$ such that for any $z\in[0,\frac{N}{n}],$ $g(z)\leq g(0)+\epsilon.$ choose $x\in(\frac{g(0)+\epsilon-h(0)}{\delta},y)$
then for any $i\in[f(n),f(n)+N]$ we have
$$[g(\frac{i-f(n)}{n})-h(0)]n^{\alpha}\leq x(i-f(n))\leq y(i-f(n))\leq [h(\frac{i-f(n)}{n})-h(0)]n^{\alpha},$$
That is because of
$$ y(i-f(n))\leq d\frac{(i-f(n))^{\nu}n^{\alpha}}{n^{\nu}},~~~~~ y(i-f(n))^{1-\nu}\leq \frac{d}{n^{\nu-\alpha}}.$$
and
$$g(\frac{i-f(n)}{n})-h(0)\leq g(0)+\epsilon-h(0)\leq x\delta< y\delta.$$
Let $x<x'<y'<y,$ we have
\beqnn &&\bfP_\mu
\left(\forall_{f(n)\leq i\leq f(n)+n}\substack{\frac{T_{i}}{n^\alpha}\in [g(\frac{i-f(n)}{n}),h(\frac{i-f(n)}{n})],\\\frac{T_{f(n)+n}}{n^\alpha}\in [a',b'],\xi_i\leq r_{n}}|T_{f(n)}=h(0)n^{\alpha}\right)
\\&=&\bfP_\mu
\Big(\forall_{f(n)\leq i\leq f(n)+n} \substack{\frac{T_{i}}{n^\alpha}\in [g(\frac{i-f(n)}{n})-h(0),h(\frac{i-f(n)}{n})-h(0)],\\ \frac{T_{f(n)+n}}{n^\alpha}+h(0)\in [a',b'],\xi_i\leq r_{n}}\Big|T_{f(n)}=0\Big)
\\&\geq&\bfP_\mu
(\forall_{f(n)\leq i\leq f(n)+N}T_{i}\in [x(i-f(n)),y(i-f(n))],\xi_i\leq r_{n}|T_{f(n)}=0)
\\&\times&\inf_{z\in[x'N,y'N]}\bfP_\mu
\Big(\forall_{N\leq i-f(n)\leq n}\substack{\frac{T_{i}}{n^\alpha}\in [g(\frac{i-f(n)}{n})-h(0),h(\frac{i-f(n)}{n})-h(0)],\\ \frac{T_{f(n)+n}}{n^\alpha}+h(0)\in [a',b'],\xi_i\leq r_{n}}\Big|T_{f(n)+N}=z\Big)
\\&\geq&\prod_{m=1}^{N}\bfP_\mu
(T_{f(n)+m}\in [x,y],\xi_{f(n)+m}\leq r_{n}|T_{f(n)+m-1}=0)
\\&\times&\inf_{z\in[x'N,y'N]}\bfP_\mu
\Big(\forall_{N\leq i-f(n)\leq n}\substack{\frac{T_{i}}{n^\alpha}\in [g(\frac{i-f(n)}{n})-h(0),h(\frac{i-f(n)}{n})-h(0)],\\ \frac{T_{f(n)+n}}{n^\alpha}\in [a'-h(0),b'-h(0)],~~~\xi_i\leq r_{n}}\Big|T_{f(n)+N}=z\Big)
\eeqnn
Let us analysis the last term of that above inequality. Notice that
$$\frac{n^\alpha}{(n-N)^\alpha}[g(\frac{i-f(n)}{n})-h(0)]=\frac{n^\alpha}{(n-N)^\alpha}[g(\frac{i-f(n)-N+N}{n-N}\frac{n-N}{n})-h(0)].$$
Hence the following two inequalities is equivalent:
$$\frac{T_{i}}{(n-N)^\alpha}\in \Big[\frac{n^\alpha}{(n-N)^\alpha}[g(\frac{i-f(n)}{n})-h(0)],\frac{n^\alpha}{(n-N)^\alpha}[h(\frac{i-f(n)}{n})-h(0)]\Big],$$
$$\small{\frac{T_{i}}{(n-N)^\alpha}\in \left[\substack{\frac{n^\alpha}{(n-N)^\alpha}\Big[g[(\frac{i-f(n)-N}{n-N}+\frac{N}{n-N})\frac{n-N}{n}]-h(0)\Big],\\\frac{n^\alpha}{(n-N)^\alpha}\Big[h[(\frac{i-f(n)-N}{n-N}+\frac{N}{n-N})\frac{n-N}{n}]-h(0)\Big]}\right]}.$$
We can see when $i$ runs from $f(n)+N$ to $f(n)+n,$ by lemma 3.2 it implies
$g[(\frac{i-f(n)-N}{n-N}+\frac{N}{n-N})\frac{n-N}{n}]$ will become $g((x+\frac{N}{n-N})\frac{n-N}{n})=g(x+(1-x)\frac{N}{n}).$ And we also can see in the sense of $L^{+\infty},$ we have
$$L^{n}_1:=\frac{n^\alpha}{(n-N)^\alpha}[g\Big((x+\frac{N}{n-N})\frac{n-N}{n}\Big)-h(0)]\stackrel{L^{+\infty}}{\longrightarrow}g(x)-h(0)$$
$$L^{n}_2:=\frac{n^\alpha}{(n-N)^\alpha}[h\Big((x+\frac{N}{n-N})\frac{n-N}{n}\Big)-h(0)]\stackrel{L^{+\infty}}{\longrightarrow}h(x)-h(0)$$
So we have
\beqnn&&\varliminf\limits_{n\rightarrow +\infty}\frac{\ln \inf_{z\in[x'N,y'N]}\bfP_\mu
\Big(\forall_{N\leq i-f(n)\leq n}\substack{\frac{T_{i}}{n^\alpha}\in [g(\frac{i-f(n)}{n})-h(0),h(\frac{i-f(n)}{n})-h(0)],\\ \frac{T_{f(n)+n}}{n^\alpha}\in [a'-h(0),b'-h(0)],~~~\xi_i\leq r_{n}}\Big|T_{f(n)+N}=z\Big)}{n^{1-2\alpha}}
\\&&~~~~~~~~~~\geq -C^{0,1}_{g,h}\sigma^2_Q\gamma(\frac{\sigma_A}{\sigma_Q}).\eeqnn
if $$\bfE([-\ln\bfP_\mu(T_{1}\in [x,y],\xi_{1}\leq A|T_{0}=0)]^8)<+\infty.$$
By 0-1 law and the fact $\lim\limits_{n\rightarrow +\infty}\frac{N}{n^{1-2\alpha}}=0, (\alpha\in(0,1/3])$ we can see \beqnn&& \frac{\ln\prod_{m=1}^{N}\bfP_\mu(T_{f(n)+m}\in [x,y],\xi_{f(n)+m}\leq r_{n}|T_{f(n)+m-1}=0)}{n^{1-2\alpha}}
\\&>&-\delta\bfE(-\ln\bfP_\mu(T_{1}\in [x,y],\xi_{1}\leq A|T_{0}=0))\eeqnn
Let $\delta\rightarrow 0,$ we complete this proof.
\section{Proof}
\noindent{\bf Proof of Theorem 2.1 (a)}
Let $M\in\bfN,$ define
$$P_n(\L):=\bfP(\forall 1\leq k\leq n,\sharp\{u\in\mathcal{T}_{M^{k}}, \forall i\leq M^{k}, V(u_i)\leq ai^{\frac{1}{3}}-\frac{K_i}{\vartheta}\}\geq v_{k-1}).$$ $z$ is a particle in this system such that $V(z)=aM^{\frac{k}{3}}-\frac{K_{M^{k}}}{\vartheta}, |z|=M^k.$ Define
$$Z_k(\L):=\sharp\Big\{u\in\mathcal{T}_{M^{(k+1)}}: \substack{\forall M^{k}<i\leq M^{(k+1)}, V(u_i)\in[(a-b)i^{\frac{1}{3}}-\frac{K_i}{\vartheta},ai^{\frac{1}{3}}-\frac{K_i}{\vartheta}],\\ \gamma(u_{i-1})\leq r_k,~~~~ u>z}\Big\},$$
$$Y_k(\L):=X_{M^k}(\L)=\sharp\{u\in\mathcal{T}_{M^{k}}, \forall i\leq M^{k}, V(u_i)\leq ai^{\frac{1}{3}}-\frac{K_i}{\vartheta}\}.$$
It is easy to see
\beqnn\frac{P_{n+1}(\L)}{P_n(\L)}&:=&\bfP_{\L}(\forall 1\leq k\leq n+1,Y_k(\L)\geq v_{k-1}|
\forall 1\leq k\leq n,Y_k(\L)\geq v_{k-1})
%\bfP_{\L}(\text{在第$e^{\lambda k}$代的$v_{k-1}+$人口中至少有$v_{k-1}$人中的一个生了符合条件的$v_k$个孩子在第$e^{\lambda (k+1)}$代})
\\&\geq&1-\bfP_{\L}(Z_n(\L)< v_n)^{\lfloor v_{n-1}\rfloor}.\eeqnn If we denote $A_k(\L):=\bfP_{\L}(Z_k(\L)\geq v_k),$  write $P_n(\L)$ as $P_n$ for simplicity, then we have
$$P_{n}\geq P_{1}\prod_{k=1}^{n-1}(1-(1-A_k)^{\lfloor v_{k-1}\rfloor})
\geq P_{1}\prod_{k=1}^{n-1}(1-e^{-A_k \lfloor v_{k-1}\rfloor}).$$
Let $v_k:=\theta \bfE_{\L}(Z_k),$ then we can see $P_1>0, {\rm \bfP_{\L}-a.s.}.$
To prove Thmorem 2.1 (a), we only need to show that
\beqlb\label{sect0-1}\sum_{i=1}^{+\infty}e^{-A_i \lfloor v_{i-1}\rfloor}<+\infty.~~{\rm \bfP_{\L}-a.s.}.\eeqlb
For simplicity, we denote $Z_n(\L)$ by $Z_n$ under the probability space $\bfP_{\L}.$
%$v_k:=\min\{\frac{(r_k)^{e^{(\lambda+1)k}-e^{k\lambda}}}{2},\theta \bfE_{\L}(Z_k)\},$

Note that $v_k:=\theta \bfE_{\L}(Z_k),$ so we can use the Paley-Zygmund inequality to get the following inequality:
~~~$$A_k:=\bfP_{\L}(Z_k\geq v_k)\geq(1-\theta)^2\frac{\bfE^2_{\L}(Z_k)}{\bfE_{\L}(Z^2_k)}.$$ Define $d_k:=M^{k+1}-M^{k},~I_i(\L):=[(a-b)i^{\frac{1}{3}}-\frac{K_i}{\vartheta},ai^{\frac{1}{3}}-\frac{K_i}{\vartheta}],$
\beqnn \bfE_{\L}(Z^2_k)&=&\bfE_{\L}\left(\sum_{\substack{u>z,v>z,\\|u|=|v|=e^{M^{k+1}}}}\1_{\{M^{k}<i\leq M^{k+1}, V(u_i)\in I_i(\L),\gamma(u_{i-1})\leq r_k\}}\right):=\sum_{j=0}^{d_k}B_{k,j}(\L)\eeqnn
%$$ \Theta_k:=\Big\{u\in \mathcal{T}_{e^{(k+1)\lambda}}: \forall M^{k}<i\leq M^{k+1}, V(u_i)\in[(a-b)i^{\frac{1}{3}}-\frac{K_i}{\vartheta},ai^{\frac{1}{3}}-\frac{K_i}{\vartheta}],\nu(u_{i-1})\leq r_{k}\Big\}.$$
By second moment method, we have
%$B_{k,j}(\L)\leq(r_k-1)\sup_{|v'|=d_k+j}\bfE_{\L}\Big(Z_k^{v'}(\Theta)\Big)\bfE_{\L}\Big(Z_k\Big)$
$B_{k,j}(\L)\leq\bfE_{\L}(Z_k)+(r_k-1)h_{k,j}(\L)\bfE_{\L}(Z_k).$
So we have
\beqlb\label{key}A_k\geq(1-\theta)^2\frac{\bfE_{\L}(Z_k)}{1+(r_k-1)h_{k,j}(\L)}.~~{\rm~\bfP_{\L}-a.s..}\eeqlb
Let
$I_{k,j}(\L):=[(a-b)(c_k+j)^{\frac{1}{3}}-\frac{K_{c_k+j}}{\vartheta},a(c_k+j)^{\frac{1}{3}}-\frac{K_{c_k+j}}{\vartheta}], c_k=M^k,$ According to the
assumption (2.2),(2.3), (2.4) and (2.5), $\{T_i\}_{i\in\bfN}$ satisfied the conditions of Lemma 3.2.  we can see%\sup_{|v'|=l_k+j}\bfE_{\L}\Big(Z_k^{v'}(\Theta)\Big)
\beqnn h_{k,j}(\L)&\leq&\sup_{x\in I_{k,j}(\L)}\bfE^{c_k+j}_{\L}\left(\sum_{|y|=d_{k+1}-j}\1_{\substack{\{\forall i\leq d_{k+1}-j,
x+V(y_i)\in[(a-b)(i+c_k+j)^{\frac{1}{3}}-\frac{K_{i+c_k+j}}{\vartheta},\\ a(i+c_k+j)^{\frac{1}{3}}-\frac{K_{i+c_k+j}}{\vartheta}]\}}}\right)
\\&=&\sup_{x\in I_{k,j}(\L)}\bfE^{c_k+j}_{\L}\left(e^{T_{d_{k+1}-j}}\1_{\substack{\{\forall i\leq d_{k+1}-j,
\vartheta x+ T_i\in[\vartheta(a-b)(i+c_k+j)^{\frac{1}{3}}-K_{c_k+j},\\ \vartheta a(i+c_k+j)^{\frac{1}{3}}-K_{c_k+j}]\}}}\right)
\\&\leq&\sup_{x\in I_{k,j}(\L)}e^{\vartheta ac_{k+1}^{\frac{1}{3}}-K_{c_k+j}-\vartheta x}\bfP^{c_k+j}_{\L}\left(\1_{\substack{\{\forall i\leq d_{k+1}-j,
\vartheta x+ T_i\in[\vartheta(a-b)(i+c_k+j)^{\frac{1}{3}}-K_{c_k+j},\\ \vartheta a(i+c_k+j)^{\frac{1}{3}}-K_{c_k+j}]\}}}\right)
\\&\leq&e^{\vartheta ac_{k+1}^{\frac{1}{3}}-\vartheta(a-b)(c_k+j)^{\frac{1}{3}}}\sup_{\substack{y\in[\vartheta(a-b)(c_k+j)^{\frac{1}{3}},\\\vartheta a(c_k+j)^{\frac{1}{3}}]}}\bfP^{c_k+j,y}_{\L}\Big(\substack{\forall i\leq d_{k+1}-j,
T_i\in[\vartheta(a-b)(i+c_k+j)^{\frac{1}{3}},\\ \vartheta a(i+c_k+j)^{\frac{1}{3}}]}\Big).\eeqnn
Where $y:=\vartheta x+K_{c_k+j}\in[\vartheta(a-b)(c_k+j)^{\frac{1}{3}},\vartheta a(c_k+j)^{\frac{1}{3}}]:=J_{k,j}.$
%$$\frac{T_i}{d^{\frac{1}{3}}_k}\in[\vartheta(a-b)(\frac{i+c_k+j}{d_k})^{\frac{1}{3}},\\ \vartheta a(\frac{i+c_k+j}{d_k})^{\frac{1}{3}}]$$
We will divide $d_k=M^{k+1}-M^{k}=M(M-1)M^{k-1}$.  Denote $K(M):=M^2-M-1,$ we have
%\beqnn &&\sum_{j=M^{k}+1}^{M^{k+1}}e^{\vartheta ac_{k+1}^{\frac{1}{3}}-\vartheta(a-b)(c_k+j)^{\frac{1}{3}}}\sup_{y\in[(a-b)(c_k+j)^{\frac{1}{3}},a(c_k+j)^{\frac{1}{3}}]}\bfP^{c_k+j,y}_{\L}\Big(\substack{\forall i\leq d_{k+1}-j,
%T_i\in[\vartheta(a-b)(i+c_k+j)^{\frac{1}{3}},\\ \vartheta a(i+c_k+j)^{\frac{1}{3}}]}\Big)
%\\&\leq&\sum_{l=0}^{M-2}e^{\vartheta ac_{k+1}^{\frac{1}{3}}-\vartheta(a-b)(c_k+j)^{\frac{1}{3}}}\sum_{j=M^k+lM^k+1}^{M^k+lM^k+M^k}\eeqnn
\beqnn &&\sum_{j=1}^{M^{k+1}-M^{k}}e^{\vartheta ac_{k+1}^{\frac{1}{3}}-\vartheta(a-b)(c_k+j)^{\frac{1}{3}}}\sup_{y\in J_{k,j}}\bfP^{c_k+j,y}_{\L}\Big(\substack{\forall i\leq d_{k+1}-j,
T_i\in[\vartheta(a-b)(i+c_k+j)^{\frac{1}{3}},\\ \vartheta a(i+c_k+j)^{\frac{1}{3}}]}\Big)
\\&\leq&\sum_{l=0}^{K(M)}e^{\vartheta ac_{k+1}^{\frac{1}{3}}-\vartheta(a-b)(c_k+lM^{k-1})^{\frac{1}{3}}}\sum_{j=lM^{k-1}+1}^{lM^{k-1}+M^{k-1}}\sup_{y\in J_{k,j}}\bfP^{c_k+j,y}_{\L}\Big(\forall_{i\leq d_{k+1}-j},
T_i\in J_{k,j+i}\Big)
\\&:=&\sum_{l=0}^{K(M)}e^{\vartheta ac_{k+1}^{\frac{1}{3}}-\vartheta(a-b)(c_k+lM^{k-1})^{\frac{1}{3}}}\sum_{j=lM^{k-1}+1}^{lM^{k-1}+M^{k-1}}H_j
\\&\leq&\sum_{l=0}^{K(M)}e^{\vartheta ac_{k+1}^{\frac{1}{3}}-\vartheta(a-b)(c_k+lM^{k-1})^{\frac{1}{3}}}M^{k-1}H_{(l+1)M^{k-1}}.\eeqnn
It is not difficult to see the $H_j$ is non-decrease by Markov property. According to corollary 3.1 and the fact that $K(M)$ is a finite fixed number (not depend on $k$), we can see
\beqnn&&\lim\limits_{k\rightarrow+\infty}\frac{\ln\sum_{l=0}^{K(M)}e^{\vartheta ac_{k+1}^{\frac{1}{3}}-\vartheta(a-b)(c_k+lM^{k-1})^{\frac{1}{3}}}M^{k-1}H_{(l+1)M^{k-1}}}{d^{1/3}_k}
\\&=&\max_{l\leq K(M)}\left[\varlimsup\limits_{k\rightarrow +\infty}\frac{\vartheta aM^{\frac{k+1}{3}}-\vartheta(a-b)(c_k+lM^{k-1})^{\frac{1}{3}}}{d^{1/3}_k}+\varlimsup\limits_{k\rightarrow +\infty}\frac{\ln H_{(l+1)M^{k-1}}}{d^{1/3}_k}\right]
\\&=&\max_{l\leq K(M)}\left[\vartheta a g_M(1)-\vartheta(a-b)g_M\Big(\frac{l}{M^2-M}\Big)-\frac{3\gamma_{\sigma}}{\vartheta^2b^2}\left(g_M(1)-g_M\left(\frac{l+1}{M^2-M}\right)\right)\right].
\eeqnn Where $g_M(x):=\Big(x+\frac{1}{M-1}\Big)^{\frac{1}{3}},~\gamma_{\sigma}:=\sigma^2_Q\gamma(\sigma_A/\sigma_Q).$ Notice that for any $l\leq K(M),$ it is true that $$g_M\left(\frac{l+1}{M^2-M}\right)-g_M\left(\frac{l}{M^2-M}\right)\leq \left(\frac{1}{M^2-M}\right)^{1/3}=\frac{g_M(0)}{M^{1/3}}.$$
Hence we have
\beqlb\label{a1u0}&&\lim\limits_{k\rightarrow+\infty}\frac{\ln\sum_{l=0}^{K(M)}e^{\vartheta ac_{k+1}^{\frac{1}{3}}-\vartheta(a-b)(c_k+lM^{k-1})^{\frac{1}{3}}}M^{k-1}H_{(l+1)M^{k-1}}}{d^{1/3}_k}\nonumber
\\&\leq& \sup_{x\in[0,1]}\left[\left(\vartheta a-\frac{3\gamma_{\sigma}}{\vartheta^2b^2}\right)g_M(1) +\left(\frac{3\gamma_{\sigma}}{\vartheta^2b^2}-\vartheta(a-b)\right)g_M\left(x\right)+\frac{g_M(0)}{M^{1/3}}\right].\eeqlb
Take $b\in(b_1, b_2),$ where $b_1,b_2$ is the two solutions of $\frac{3\gamma_{\sigma}}{\vartheta^2b^2}-\vartheta(a-b)=0.$ Hence we have
\beqnn
&&\sup_{x\in[0,1]}\left[\left(\vartheta a-\frac{3\gamma_{\sigma}}{\vartheta^2b^2}\right)g_M(1) +\left(\frac{3\gamma_{\sigma}}{\vartheta^2b^2}-\vartheta(a-b)\right)g_M\left(x\right)+\frac{g_M(0)}{M^{1/3}}\right]
\\&=&\left(\vartheta a-\frac{3\gamma_{\sigma}}{\vartheta^2b^2}\right)g_M(1) +\left(\frac{3\gamma_{\sigma}}{\vartheta^2b^2}-\vartheta(a-b)\right)g_M\left(0\right)+\frac{g_M(0)}{M^{1/3}}
\\&=&\left[\left(\vartheta a-\frac{3\gamma_{\sigma}}{\vartheta^2b^2}\right)M^{1/3} +\left(\frac{3\gamma_{\sigma}}{\vartheta^2b^2}+\vartheta b-\vartheta a\right)+\frac{1}{M^{1/3}}\right]g_M(0).
\eeqnn
From the above we can see $$\sup_{x\in[0,1]}\left[\left(\vartheta a-\frac{3\gamma_{\sigma}}{\vartheta^2b^2}\right)g_M(1) +\left(\frac{3\gamma_{\sigma}}{\vartheta^2b^2}-\vartheta(a-b)\right)g_M\left(x\right)+\frac{g_M(0)}{M^{1/3}}\right]>0.$$ when $M$ is large enough. That is to say, if we take $r_n=e^{n^{1/4}},$ then we have
\beqlb\label{a1u}&&\lim\limits_{k\rightarrow+\infty}\frac{\ln[1+(r_k-1)h_{k,j}(\L)]}{d^{1/3}_k}\nonumber
\\&\leq&\sup_{x\in[0,1]}\left[\left(\vartheta a-\frac{3\gamma_{\sigma}}{\vartheta^2b^2}\right)g_M(1) +\left(\frac{3\gamma_{\sigma}}{\vartheta^2b^2}-\vartheta(a-b)\right)g_M\left(x\right)+\frac{g_M(0)}{M^{1/3}}\right]\eeqlb

%$$\vartheta a(\frac{M}{M-1})^{\frac{1}{3}}-\vartheta(a-b)\sqrt[3]{\frac{l+1}{M-1}}$$
%$$\vartheta a(\frac{M}{M-1})^{\frac{1}{3}}-\vartheta(a-b)\sqrt[3]{\frac{l+1}{M-1}}--\frac{3\gamma(\beta)\sigma^2}{b^2\vartheta^2}(\sqrt[3]{\frac{M}{M-1}}-\sqrt[3]{\frac{l+1}{M-1}})$$
%so let $a_M=\lfloor M^{2/5}\rfloor$
Let us turn to the lower bound of the $\bfE_{\L}(Z_k).$ We write $\bfE^{t,k}_{\L}(\cdot):=\bfE_{\L}(\cdot|V(z)=k,|z|=t).$ Then we have
\beqnn \bfE_{\L}(Z_k)&=&\bfE^{c_k,ac_k^{\frac{1}{3}}-\frac{K_{c_k}}{\vartheta}}_{\L}\left(\sum_{u>z,|u|=e^{\lambda(k+1)}}\1_{\{e^{k\lambda}<i\leq e^{(k+1)\lambda}, V(u_i)\in I_i(\L),\gamma(u_{i-1})\leq r_k\}}\right)
\\&=&\bfE^{c_k,\vartheta ac_k^{\frac{1}{3}}-K_{c_k}}_{\L}\left(e^{T_{d_k}-[\vartheta ac_k^{\frac{1}{3}}-K_{c_k}]}\1_{\substack{0<i\leq d_k,\xi_i\leq r_k,\\ T_i\in[\vartheta(a-b)(i+c_k)^{\frac{1}{3}}-K_{c_k},\vartheta a(i+c_k)^{\frac{1}{3}}-K_{c_k}]}}\right)
\\&=&\bfE^{c_k,\vartheta ac_k^{\frac{1}{3}}}_{\L}(e^{T_{d_k}-\vartheta ac_k^{\frac{1}{3}}}\1_{0<i\leq d_k,\xi_i\leq r_k,T_i\in[\vartheta(a-b)(i+c_k)^{\frac{1}{3}},\vartheta a(i+c_k)^{\frac{1}{3}}]})
\\&\geq&e^{(\vartheta a-\epsilon)c_{k+1}^{\frac{1}{3}}-\vartheta ac_k^{\frac{1}{3}}}\bfP^{c_k,\vartheta ac_k^{\frac{1}{3}}}_{\L}\left(0<i\leq d_k,\xi_i\leq r_k,\substack{T_i\in[\vartheta(a-b)(i+c_k)^{\frac{1}{3}},\vartheta a(i+c_k)^{\frac{1}{3}}]\\,
T_{d_k}\in[\vartheta(a-\epsilon)c_{k+1}^{\frac{1}{3}},\vartheta ac_{k+1}^{\frac{1}{3}}]}\right).\eeqnn
%\beqnn&&\bfP^{c_k,\vartheta ac_k^{\frac{1}{3}}}_{\L}(0<i\leq d_k,\xi_i\leq r_k,\substack{T_i\in[\vartheta(a-b)(i+c_k)^{\frac{1}{3}},\vartheta a(i+c_k)^{\frac{1}{3}}]\\,
%T_{d_k}\in[\vartheta(a-\epsilon)c_{k+1}^{\frac{1}{3}},\vartheta ac_{k+1}^{\frac{1}{3}}]})
%\\&\geq&\bfP^{c_k,\vartheta ac_k^{\frac{1}{3}}}_{\L}(i\leq \varepsilon \sqrt[3]{d_k},\xi_i\leq r_k,T_i\in[\vartheta ac_{k}^{\frac{1}{3}}-v_2i, \vartheta ac_{k}^{\frac{1}{3}}-v_1i])\eeqnn
By condition 4, we can utilize the Corollary 3.2 to get the following limit.
%$$ \Theta_k:=\Big\{u\in \mathcal{T}_{e^{(k+1)\lambda}}: \forall e^{k\lambda}<i\leq e^{(k+1)\lambda}, V(u_i)\in[(a-b)i^{\frac{1}{3}}-\frac{K_i}{\vartheta},ai^{\frac{1}{3}}-\frac{K_i}{\vartheta}],\nu(u_{i-1})\leq r_{k}\Big\}.$$
%$$\vartheta(a-b)\sqrt[3]{i+c_k}\leq\vartheta a\sqrt[3]{c_k}-v_2 i\leq \vartheta a\sqrt[3]{c_k}+v_1 i\leq\vartheta(a-b)\sqrt[3]{c_k}$$
\beqnn&&\varliminf\limits_{k\rightarrow+\infty}\frac{\ln \bfE_{\L}(Z_k)}{d_k}
\\&\geq&(a\vartheta-\epsilon)(\frac{M}{M-1})^{\frac{1}{3}}-a\vartheta(\frac{1}{M-1})^{\frac{1}{3}}-\sigma^2\gamma(\beta)\int_{0}^1[b\vartheta(x+\frac{1}{M-1})^{\frac{1}{3}}]^{-2}dx
\\&=&a\vartheta(g_{M}(1)-g_{M}(0))-\frac{3\gamma(\beta)\sigma^2}{b^2\vartheta^2}(g_{M}(1)-g_{M}(0))>0.\eeqnn
Therefore
$$\varliminf\limits_{k\rightarrow+\infty}\frac{\ln \lfloor\bfE_{\L}(Z_k)\rfloor}{d_k}\geq a\vartheta (g_{M}(1)-g_{M}(0))-\frac{3\gamma(\beta)\sigma^2}{b^2\vartheta^2}(g_{M}(1)-g_{M}(0)).$$
Combining with \eqref{a1u} we can see
\beqnn&&\lim\limits_{k\rightarrow+\infty}\frac{\ln A_k \lfloor v_{k-1}\rfloor}{d_k}
\\&\geq& f(0)(\vartheta a-\vartheta b-\frac{3\gamma(\beta)\sigma^2}{b^2\vartheta^2})-f(0)o(M). \eeqnn
We can see the when $b=\frac{\sqrt[3]{6\gamma(\beta)\sigma^2}}{\vartheta},$ $\vartheta b+\frac{3\gamma(\beta)\sigma^2}{b^2\vartheta^2}$ take its minimum value $\frac{3}{2}\sqrt[3]{6\gamma(\beta)\sigma^2}.$ So if $a>\frac{3\sqrt[3]{6\gamma(\beta)\sigma^2}}{2\vartheta},$ $\lim\limits_{k\rightarrow+\infty}\frac{\ln A_k \lfloor v_{k-1}\rfloor}{d_k}>0.~{\rm~\bfP_{\L}-a.s..}$ Then we can see (4.1) holds, thus we complete the proof of Theorem 2.1 (a).

\bigskip

――――――――――――――――――――――――――――――――――――――――――――――――――――――――――――――――――――――――――――――――――――――

\bigskip

\noindent{\bf Proof of Theorem 2.1 (b): the upper bound}

Let $g(i)=ai^{\frac{1}{3}},$ $f(\cdot): [0,1]\rightarrow[0,+\infty)$ is a continue non-negative function. Define
 $$Z_n:=\sharp\{|u|=n, \forall i\leq n, V(u_i)\leq ai^{\frac{1}{3}}-\frac{K_i}{\vartheta}\}$$
$$\bfP_{\L}(Z_n>0)=\bfP_{\L}(\exists u: |u|=n, \forall i\leq n, V(u_i)\leq ai^{\frac{1}{3}}-\frac{K_i}{\vartheta})\leq \sum_{j=1}^{n}H_j+H.$$
Where $$H_j:=\bfP_{\L}\Big(\exists |u|=j: \forall i<j, \substack{V(u_i)\in[ai^{\frac{1}{3}}-n^{\frac{1}{3}}f(\frac{i}{n})-\frac{K_i}{\vartheta},g(i)-\frac{K_i}{\vartheta}],\\ V(u_j)\leq aj^{\frac{1}{3}}-n^{\frac{1}{3}}f(\frac{j}{n})-\frac{K_j}{\vartheta}}\Big).$$
$$H:=\bfP_{\L}\Big(\exists |u|=n: \forall i\leq n, V(u_i)\in[ai^{\frac{1}{3}}-n^{\frac{1}{3}}f(\frac{i}{n})-\frac{K_i}{\vartheta},g(i)-\frac{K_i}{\vartheta}]\Big).$$
By Markov inequality and many to one formula, we have
\beqnn H_j&\leq& \bfE_{\L}\Big(\sum_{|u|=j}\1_{\{\forall i<j, V(u_i)\in[ai^{\frac{1}{3}}-n^{\frac{1}{3}}f(\frac{i}{n})-\frac{K_i}{\vartheta},g(i)-\frac{K_i}{\vartheta}],V(u_j)\leq aj^{\frac{1}{3}}-n^{\frac{1}{3}}f(\frac{j}{n})-\frac{K_j}{\vartheta}\}}\Big)
\\&=&\bfE_{\L}\Big(e^{T_j}\1_{\{\forall i<j, S_i\in[ai^{\frac{1}{3}}-n^{\frac{1}{3}}f(\frac{i}{n})-\frac{K_i}{\vartheta},g(i)-\frac{K_i}{\vartheta}],S_j\leq aj^{\frac{1}{3}}-n^{\frac{1}{3}}f(\frac{j}{n})-\frac{K_j}{\vartheta}\}}\Big)
\\&\leq& e^{\vartheta aj^{\frac{1}{3}}-\vartheta n^{\frac{1}{3}}f(\frac{j}{n})}\bfP_{\L}\Big(\forall i<j, S_i\in[ai^{\frac{1}{3}}-n^{\frac{1}{3}}f(\frac{i}{n})-\frac{K_i}{\vartheta},g(i)-\frac{K_i}{\vartheta}]\Big)
\\&=& e^{\vartheta aj^{\frac{1}{3}}-\vartheta n^{\frac{1}{3}}f(\frac{j}{n})}\bfP_{\L}\Big(\forall i<j, T_i\in[\vartheta ai^{\frac{1}{3}}-\vartheta n^{\frac{1}{3}}f(\frac{i}{n}),\vartheta ai^{\frac{1}{3}}]\Big).\eeqnn
 $$H\leq e^{\vartheta g(n)}\bfP_{\L}\Big(\forall i\leq n, S_i\in[ai^{\frac{1}{3}}-n^{\frac{1}{3}}f(\frac{i}{n})-\frac{K_i}{\vartheta},g(i)-\frac{K_i}{\vartheta}]\Big).$$
 For the monotonicity of $\bfP_{\L}(Z_n>0)$ we only need to consider $n:=Nk, l\in[0, N-1]\cap\bfN.$ $f_{l,N}=\inf_{x\in[\frac{l}{N},\frac{l+1}{N}]}f(x).$
 $$\sum_{i=kl+1}^{k(l+1)}H_i\leq ke^{\vartheta a((l+1)k)^{\frac{1}{3}}-\vartheta n^{\frac{1}{3}}f_{l,N}}\bfP_{\L}\big(\forall i\leq lk, T_i\in[\vartheta ai^{\frac{1}{3}}-\vartheta n^{\frac{1}{3}}f(\frac{i}{n})  , \vartheta ai^{\frac{1}{3}}  ]\big)$$
 By Lemma 3.2, we have
 \beqnn\lim\limits_{k\rightarrow+\infty}\frac{\ln(\sum\limits_{i=kl+1}^{k(l+1)}H_i)}{n^{\frac{1}{3}}}&\leq&\vartheta a(\frac{l+1}{N})^{\frac{1}{3}}
 -\vartheta f_{l,N}-\gamma(\beta)\sigma^2\int_0^{\frac{l}{N}}(\vartheta f(x))^{-2}dx
 \\&\leq&\vartheta a(\frac{l}{N})^{\frac{1}{3}}
 -\vartheta f(\frac{l}{N})-\gamma(\beta)\sigma^2\int_0^{\frac{l}{N}}(\vartheta f(x))^{-2}dx+2\varepsilon
 \eeqnn
 \beqnn\lim\limits_{k\rightarrow+\infty}\frac{\ln H}{n^{\frac{1}{3}}}\leq\vartheta a -\gamma(\beta)\sigma^2\int_0^{1}(\vartheta f(x))^{-2}dx\eeqnn
 In conclusion, We have $$\varlimsup\limits_{n\rightarrow+\infty}\frac{\ln\bfP_{\L}(Z_n>0)}{n^\frac{1}{3}}\leq \sup_{0\leq\alpha\leq 1} [\vartheta a\alpha^{\frac{1}{3}}
 -\vartheta f(\alpha)-\gamma(\beta)\sigma^2\int_0^{\alpha}(\vartheta f(x))^{-2}dx].$$

 \noindent{\bf Proof of Theorem 2.1 (b): the lower bound}
The method of this lower bound is similar with Theorem 2.1 (a). we can get
$$\varliminf\limits_{n\rightarrow+\infty}\frac{\ln\bfP_{\L}(Z_n>0)}{n^\frac{1}{3}}\geq \inf_{0\leq\alpha\leq 1} [\vartheta a\alpha^{\frac{1}{3}}
 -\vartheta f(\alpha)-\gamma(\beta)\sigma^2\int_0^{\alpha}(\vartheta f(x))^{-2}dx].$$
According to the discussion of Proposition 3.2-Proposition 3.6 in \cite{BJ2012}, we know there exists a function $f(x)$ define on $[0,1]$ such that $\int_0^1f^{-2}(x)dx<+\infty, f(1)=0, f(0)>0$ and $$\vartheta a\alpha^{1/3}-\vartheta f(\alpha)-\frac{\gamma(\beta)\sigma^2}{\vartheta^2}\int_{0}^{\alpha}(f(x))^{-2}dx= -\vartheta f(0), \forall \alpha\in[0,1]. $$
Combining with (4.2) and (4.3), we complete the proof of Theorem 2.1 (b).

%\beqnn&&\lim\limits_{k\rightarrow+\infty}\frac{\ln\bfE_{\L}(Z_n(\Theta))}{n^\frac{1}{3}}
%\\&\geq&\vartheta a-\vartheta\varepsilon-\vartheta f(1)-\frac{\gamma(\beta)\sigma^2}{\vartheta^2}\int_{0}^{1}(\varepsilon+ f(x))^{-2}dx
%\eeqnn

\noindent{\bf Proof of Theorem 2.2: the upper bound}

Let $f(n,k):=-d\sqrt[3]{n-k}.$ We know
$$\bfP_{\L}(Z_n>0)=\bfP_{\L}(\exists u: |u|=n, V(u_i)\leq g(i)-\frac{K_i}{\vartheta})\leq \sum_{j=1}^{n}H_j+H.$$
Where $$H_j:=\bfP_{\L}\Big(\exists |u|=j: \forall i<j, V(u_i)\in[f(n,i)-\frac{K_i}{\vartheta},g(i)-\frac{K_i}{\vartheta}],V(u_j)\leq f(n,j)-\frac{K_j}{\vartheta}\Big).$$
$$H:=\bfP_{\L}\Big(\exists |u|=n: \forall i\leq n, V(u_i)\in[f(n,i)-\frac{K_i}{\vartheta},g(i)-\frac{K_i}{\vartheta}]\Big).$$
By Markov inequality and many to one formula, we have
\beqnn H_j&\leq& \bfE_{\L}\Big(\sum_{|u|=j}\1_{\{\forall i<j, V(u_i)\in[f(n,i)-\frac{K_i}{\vartheta},g(i)-\frac{K_i}{\vartheta}],V(u_j)\leq f(n,j)-\frac{K_j}{\vartheta}\}}\Big)
\\&=&\bfE_{\L}\Big(e^{T_j}\1_{\{\forall i<j, S_i\in[f(n,i)-\frac{K_i}{\vartheta},g(i)-\frac{K_i}{\vartheta}],S_j\leq f(n,j)-\frac{K_j}{\vartheta}\}}\Big)
\\&\leq& e^{\vartheta f(n,j)}\bfP_{\L}\Big(\forall i<j, S_i\in[f(n,i)-\frac{K_i}{\vartheta},g(i)-\frac{K_i}{\vartheta}]\Big).\eeqnn
 $$H\leq e^{\vartheta g(n)}\bfP_{\L}\Big(\forall i\leq n, S_i\in[f(n,i)-\frac{K_i}{\vartheta},g(i)-\frac{K_i}{\vartheta}]\Big).$$
For the monotonicity of $\bfP_{\L}(Z_n>0)$ we only need to consider $n:=Nk, l\in[0, N-1]\cap\bfN.$
 $$\sum_{i=kl+1}^{k(l+1)}H_i\leq ke^{-\vartheta d\sqrt[3]{Nk-k(l+1)}}\bfP_{\L}\big(\forall i\leq lk, T_i\in[\vartheta f(n,i), \vartheta g(i)]\big)$$
 By the random version of Mogul'ski\v{\i} estimation, we have
 \beqnn\lim\limits_{k\rightarrow+\infty}\frac{\ln(\sum\limits_{i=kl+1}^{k(l+1)}H_i)}{n^{\frac{1}{3}}}&\leq& -\vartheta d\sqrt[3]{1-\frac{l+1}{N}}
 -\gamma(\beta)\sigma^2\int_0^{\frac{l}{N}}(\vartheta d\sqrt[3]{1-x})^{-2}dx
 \\&\leq& -\vartheta d\Big(1-\frac{l+1}{N}\Big)^\frac{1}{3}-\frac{3\gamma(\beta)\sigma^2}{\vartheta^2d^2}\Big[1-\Big(1-\frac{l}{N}\Big)^{\frac{1}{3}}\Big].\eeqnn
 \beqnn\lim\limits_{k\rightarrow+\infty}\frac{\ln H}{n^{\frac{1}{3}}}\leq -\frac{3\gamma(\beta)\sigma^2}{\vartheta^2d^2}.\eeqnn
 Let $d:=\frac{\sqrt[3]{3\gamma(\beta)\sigma^2}}{\vartheta},$ we have $\limsup\limits_{n\rightarrow+\infty}\frac{\ln\bfP_{\L}(Z_n>0)}{n^\frac{1}{3}}\leq -\sqrt[3]{3\gamma(\beta)\sigma^2}.$

 \noindent{\bf Proof of Theorem 2.2: the lower bound}

Define $$ \Theta:=\Big\{u\in \mathcal{T}: \forall 1\leq i\leq |u|, V(u_i)\in[f(|u|,i)-\varepsilon |u|^{\frac{1}{3}}-\frac{K_i}{\vartheta},g(i)-\frac{K_i}{\vartheta}],\gamma(u_{i-1})\leq r_{|u|}
\Big\}.$$
$Z_n(\Theta):=\sum_{|u|=n}\1_{\{u\in\Theta\}}.$ $\bfP_{\L}(Z_n>0)\geq \bfP_{\L}(Z_n(\Theta)>0)\geq \frac{\bfE_{\L}(Z_n(\Theta))^2}{\bfE_{\L}(Z^2_n(\Theta))}.$

In this proof we write the so-called second moment method in more detail.
\beqnn\sum\limits_{|u|=n}\sum\limits_{k=0}^{n-1}\1_{\{u\in\Theta\}}Z_n^{u_k}(\Theta,u_{k+1})&=&
\sum\limits_{k=0}^{n-1}\sum\limits_{|u|=n}\1_{\{u\in\Theta\}}Z_n^{u_k}(\Theta,u_{k+1})\\&=&
\sum\limits_{k=0}^{n-1}\sum\limits_{|u|=k+1}Z_n^{u_{k+1}}(\Theta)Z_n^{u_k}(\Theta,u_{k+1}),\eeqnn
$$\sum\limits_{|u|=n}\1_{\{u\in\Theta\}}Z_n^{u_k}(\Theta,u_{k+1})=
\sum\limits_{|u|=k+1}Z_n^{u_{k+1}}(\Theta)Z_n^{u_k}(\Theta,u_{k+1}).$$
On the other hand
$$\sum\limits_{|u|=n}\sum\limits_{k=0}^{n-1}\1_{\{u\in\Theta\}}Z_n^{u_k}(\Theta,u_{k+1})=\sum\limits_{|u|=n}\1_{\{u\in\Theta\}}(Z_n(\Theta)-\1_{\{u\in\Theta\}})
=Z^2_n(\Theta)-Z_n(\Theta).$$
So \beqnn Z^2_n(\Theta)&=&Z_n(\Theta)+\sum\limits_{k=0}^{n-1}\sum\limits_{|u|=k+1}Z_n^{u_{k+1}}(\Theta)Z_n^{u_k}(\Theta,u_{k+1})
\\&=&Z_n(\Theta)+\sum\limits_{k=1}^{n}\sum\limits_{|v|=k}Z_n^{v}(\Theta)Z_n^{\overleftarrow{v}}(\Theta,v).\eeqnn
\beqnn\bfE_{\L}(\sum\limits_{|v|=k}Z_n^{v}(\Theta)Z_n^{\overleftarrow{v}}(\Theta,v))&=&\bfE_{\L}\Big(\bfE_{\L}\big(\sum\limits_{|v|=k}Z_n^{v}(\Theta)Z_n^{\overleftarrow{v}}(\Theta,v)|\mathcal{F}_k\big)\Big)
\\&\leq&\bfE_{\L}\Big(\sum\limits_{|v|=k}\bfE_{\L}\Big[\sum_{v'=bro(v)}Z_n^{v}(\Theta)Z_n^{v'}(\Theta)|\mathcal{F}_k\Big]\Big)
\\&\leq&\bfE_{\L}\Big(\sum\limits_{|v|=k}\sum_{v'=bro(v)}\bfE_{\L}\Big[Z_n^{v}(\Theta)Z_n^{v'}(\Theta)|\mathcal{F}_k\Big]\Big)
\\&\leq&\bfE_{\L}\Big(\sum\limits_{|v|=k}\sum_{v'=bro(v)}\bfE_{\L}\Big[Z_n^{v}(\Theta)|\mathcal{F}_k\Big]\bfE_{\L}\Big[Z_n^{v'}(\Theta)|\mathcal{F}_k\Big]\Big)
\\&\leq&\bfE_{\L}\Big(\sum\limits_{|v|=k}\bfE_{\L}\Big[Z_n^{v}(\Theta)|\mathcal{F}_k\Big]\sum_{v'=bro(v)}\bfE_{\L}\Big[Z_n^{v'}(\Theta)|\mathcal{F}_k\Big]\Big)
\\&\leq&\bfE_{\L}\Big((r_n-1)\sup_{|v'|=k}\bfE_{\L}\Big[Z_n^{v'}(\Theta)|\mathcal{F}_k\Big]\sum\limits_{|v|=k}\bfE_{\L}\Big[Z_n^{v}(\Theta)|\mathcal{F}_k\Big]\Big)
\\&\leq&(r_n-1)\sup_{|v'|=k}\bfE_{\L}\Big(Z_n^{v'}(\Theta)\Big)\bfE_{\L}\Big(\sum\limits_{|v|=k}\bfE_{\L}\Big[Z_n^{v}(\Theta)|\mathcal{F}_k\Big]\Big)
\\&\leq&(r_n-1)\sup_{|v'|=k}\bfE_{\L}\Big(Z_n^{v'}(\Theta)\Big)\bfE_{\L}\Big(\sum\limits_{|v|=k}Z_n^{v}(\Theta)\Big)\eeqnn
$$\bfE_{\L}(Z^2_n(\Theta))=\bfE_{\L}\left(Z_n(\Theta)\right)\Big[1+(r_n-1)\sum\limits_{k=1}^{n}\sup_{|v'|=k}\bfE_{\L}\Big(Z_n^{v'}(\Theta)\Big)\Big]$$
Let $I_{n,k}(\L)=[f(n,k)-\varepsilon n^{\frac{1}{3}}-\frac{K_k}{\vartheta},g(k)-\frac{K_k}{\vartheta}],$ then we have
\beqnn&&\sup_{|v'|=k}\bfE_{\L}\Big(Z_n^{v'}(\Theta)\Big)
\\&\leq&\sup_{V(v')\in I_{n,k}(\L)}\bfE_{\L}\Big(\sum_{|u|=n,u_k=v'}\1_{\{u\in\Theta\}}\Big)
\\&\leq&\sup_{V(v')\in I_{n,k}(\L)}\bfE_{\L}\Big(\sum_{|u|=n,u_k=v'}\1_{\{\forall i\leq n-k,V(u_{k+i})\in[f(n,i+k)-\frac{K_{i+k}}{\vartheta},g(i+k)-\frac{K_{i+k}}{\vartheta}]\}}\Big)
\\&=&\sup_{x\in I_{n,k}(\L)}\bfE^k_{\L}\Big(\sum_{|y|=n-k}\1_{\{\forall i\leq n-k,V(y_i)\in[-x+f(n,i+k)-\frac{K_{i+k}}{\vartheta},-x+g(i+k)-\frac{K_{i+k}}{\vartheta}]\}}\Big)
\\&=&\sup_{x\in I_{n,k}(\L)}\bfE^k_{\L}\Big(e^{T_{n-k}}\1_{\{\forall i\leq n-k,\frac{T_i}{\vartheta}\in[-x+f(n,i+k)-\frac{K_{k}}{\vartheta},-x+g(i+k)-\frac{K_{k}}{\vartheta}]\}}\Big).\eeqnn
Furthermore, we have
\beqnn&&\sup_{|v'|=k}\bfE_{\L}\Big(Z_n^{v'}(\Theta)\Big)
\\&\leq& e^{\vartheta g(n)-\vartheta f(n,k)}\sup_{x\in I_{n,k}(\L)}\bfP^k_{\L}\Big(\forall i\leq n-k,\frac{T_i}{\vartheta}+x+\frac{K_{k}}{\vartheta}\in[f(n,i+k),g(i+k)]\Big)
\\&\leq&e^{\vartheta g(n)-\vartheta f(n,k)}\sup_{x\in [\vartheta f(n,k),\vartheta g(k)]}\bfP^{k}_{\L}\Big(\forall i\leq n-k,T_i\in
[\vartheta f(n,i+k),\vartheta g(i+k)]|T_0=x\Big)
\\&\leq&e^{-\vartheta f(n,k)} \sup_{x\in [\vartheta f(n,k),0]}\bfP^{k}_{\L}\Big(\forall i\leq n-k,T_i\in
[\vartheta f(n,i+k),0]|T_0=x\Big).
\eeqnn
Recalling that under $\bfP^{k,x}_{\L}$ we have agreed that $T_n=\vartheta S_n+\sum\limits_{i=k+1}^{k+n}\kappa_i(\vartheta).$  Then we have
\beqnn&&\bfP_{\L}(Z_n>0)
\\&\geq&\frac{\bfE_{\L}(Z_n(\Theta))}{1+(r_n-1)\sum_{j=1}^{n}e^{-\vartheta f(n,j)}\sup_{x\in [\vartheta f(n,j),0]}\bfP^{j,x}_{\L}\Big(\forall i\leq n-j,T_i\in
[\vartheta f(n,i+j),0]\Big)}.\eeqnn
Similar with the corresponding part of the proof of Theorem 2.1(a),
 take $\vartheta d=\frac{3\gamma(\beta)\sigma^2}{d^2\vartheta^2},$ we can get $\liminf\limits_{n\rightarrow\infty}\frac{\bfP_{\L}(Z_n>0)}{n^\frac{1}{3}}\geq-\sqrt[3]{3\gamma(\beta)\sigma^2}.$ This is the end of the proof of Theorem 2.2.

――――――――――――――――――――――――――――――――――――――――――――――――――――――――――――

 \ack
I want to thank my supervisor Wenming Hong for his constant concern on my work and giving me a good learning environment. I also want to thank Bastien Mallein
for his great help on \cite{LY201801}, which is a basis of this work.

\end{document}